\documentclass[12pt]{amsart}
\usepackage{amssymb}
\usepackage[all]{xy}
\usepackage{lscape}

\theoremstyle{plain} 
\newtheorem{theor}[equation]{Theorem}
\newtheorem{cor}[equation]{Corollary}
\newtheorem{lem}[equation]{Lemma}
\newtheorem{conjecture}[equation]{Conjecture}
\newtheorem{proposition}[equation]{Proposition}

\theoremstyle{definition}
\newtheorem{defin}[equation]{Definition}

\theoremstyle{remark}

\newtheorem{ex}[equation]{Example}

\newcommand{\BV}{\Delta}
\newcommand{\sectiontriviale}{s}
\def\build#1_#2^#3{\mathrel{\mathop{\kern0pt#1}\limits_{#2}^{#3}}}

\begin{document}
\begin{abstract}
Let $M$ be a compact oriented $d$-dimensional smooth manifold.
Chas and Sullivan have defined a structure of Batalin-Vilkovisky
algebra on $\mathbb{H}_*(LM)$.
Extending work of Cohen, Jones and Yan,
we compute this Batalin-Vilkovisky algebra structure when
$M$ is a sphere $S^d$, $d\geq 1$.
In particular, we show that $\mathbb{H}_*(LS^2;\mathbb{F}_2)$ and the Hochschild cohomology
$HH^{*}(H^*(S^2);H^*(S^2))$ are surprisingly not isomorphic as Batalin-Vilkovisky algebras,
although we prove that, as expected, the underlying Gerstenhaber
algebras are isomorphic.
The proof requires the knowledge of the Batalin-Vilkovisky algebra
$H_*(\Omega^2 S^3;\mathbb{F}_2)$ that we compute in the Appendix.
\end{abstract}
\title{\bf String topology for spheres.}
\author{Luc Menichi*\\
With an appendix by Gerald Gaudens and Luc Menichi}

\address{UMR 6093 associ\'ee au CNRS\\
Universit\'e d'Angers, Facult\'e des Sciences\\
2 Boulevard Lavoisier\\49045 Angers, FRANCE}
\email{firstname.lastname at univ-angers.fr}
\thanks{*The author was partially supported by the Mathematics
  Research Center of Stanford University}\maketitle
\begin{center}
\textit{Dedicated to Jean-Claude Thomas, on the occasion of his
60th birthday}
\end{center}
\section{Introduction}
Let $M$ be a compact oriented $d$-dimensional smooth manifold.
Denote by $LM:=map(S^1,M)$ the free loop space on $M$.
In 1999, Chas and Sullivan~\cite{Chas-Sullivan:stringtop} have shown that the shifted free loop homology
$\mathbb{H}_*(LM):=H_{*+d}(LM)$ has a structure of Batalin-Vilkovisky algebra
(Definition~\ref{definition algebre de Batalin-Vilkovisky}). In particular, they showed
that $\mathbb{H}_*(LM)$ is a Gerstenhaber algebra (Definition~\ref{definition algebre de Gerstenhaber}).
This Batalin-Vilkovisky algebra has been computed when $M$ is
a complex Stiefel manifold~\cite{Tamanoi:BVLalg} and very recently
over $\mathbb{Q}$ when $M$ is a $K(\pi,1)$~\cite{Vaintrob:stBVaHcohGoldb}.
In this paper, we compute the Batalin-Vilkovisky algebra $\mathbb{H}_*(LM;{\Bbbk})$ when $M$ is a sphere $S^n$, $n\geq 1$
over any commutative ring ${\Bbbk}$ (Theorems~\ref{string topology du cercle}, \ref{string topology des spheres impairs},
\ref{string topology spheres paires sur les entiers}, \ref{homologie des lacets sur la sphere modulo 2} and
\ref{homologie des lacets sur la sphere modulo Z}).

In fact, few  calculations of this Batalin-Vilkovisky algebra
structure
or even of the underlying Gerstenhaber
algebra structure 
have been done because the following conjecture
has not yet been proved.
\begin{conjecture}(due to ~\cite[``dictionary'' p. 5]{Chas-Sullivan:stringtop} or~\cite{CohJon:ahtrostringtopology}?)\label{conjecture}

\noindent If $M$ is simply connected then there is an isomorphism of Gerstenhaber algebras
$\mathbb{H}_*(LM)\cong HH^*(S^*(M);S^*(M))$ between the free loop space homology and the Hochschild cohomology of the
algebra of singular cochains on $M$.
\end{conjecture}
In~\cite{CohJon:ahtrostringtopology,Cohen:multpropAtiyahduality},
Cohen and Jones proved that
there is an isomorphism of graded algebras over any field
$$\mathbb{H}_*(LM)\cong HH^*(S^*(M);S^*(M)).$$
Over the reals or over the rationals, two proofs of this isomorphism of graded algebras
have been given by Merkulov~\cite{Merkulov:DeRhammfst} and F\'elix, Thomas, Vigu\'e-Poirrier~\cite{Felix-Vigue-Thomas:ratstringtop}.
Motivated by this conjecture, Westerland~\cite{Westerland:DyerLoitstosaps} has computed the Gerstenhaber algebra
$
HH^*(S^*(M;\mathbb{F}_2);S^*(M;\mathbb{F}_2))
$
when $M$ is a sphere or a projective space.

What about the Batalin-Vilkovisky algebra structure?

Suppose that $M$ is formal over a field, then since the Gerstenhaber algebra structure on
Hochschild cohomology is preserves by quasi-isomorphism of
algebras~\cite[Theorem 3]{Felix-Menichi-Thomas:GerstduaiHochcoh},
we obtain an isomorphism of Gerstenhaber algebras

\begin{equation}\label{invariance cohomologie de Hochschild}
HH^*(S^*(M);S^*(M))\cong HH^*(H^*(M);H^*(M)).
\end{equation}
Poincar\'e duality induces an isomorphism of $H^*(M)$-modules
$$\Theta:H^*(M)\rightarrow H^*(M)^{\vee}.$$
Therefore, we obtain the isomorphism
$$HH^*(H^*(M);H^*(M))\cong HH^*(H^*(M);H^*(M)^{\vee})$$
and the Gerstenhaber algebra structure on $HH^*(H^*(M);H^*(M))$ extends
to a Batalin-Vilkovisky algebra~\cite{Tradler:bvalgcohiiip,
  MenichiL:BValgaccoHa,Kaufmann:procvDelignecvc} (See above Proposition~\ref{Cohomologie de Hochschild des spheres modulo 2} for details).
This Batalin-Vilkovisky algebra structure is further extended
in~\cite{Tradler-Zeinalian:ontcyclicDelignec,Costello:tcftcalabiyaucat,Kaufmann:modsahcfaco,Kontsevich-Soibelman:noainfinityalg}
to a richer algebraic structure.
It is natural to conjecture that this Batalin-Vilkovisky algebra on
$HH^*(H^*(M);H^*(M))$ is isomorphic to the Batalin-Vilkovisky algebra $\mathbb{H}_*(LM)$.
We show (Corollary~\ref{homologie des lacets different de cohomologie de Hochschild comme BV-algebre})
that this is not the case over $\mathbb{F}_2$ when $M$ is the sphere $S^2$.
See~\cite[Comments 2. Chap. 1]{Cohen-Hess-Voronov:stringtopacyclhom} or
the papers of Tradler and
Zeinalian~\cite{Tradler:bvalgcohiiip,Tradler-Zeinalian:ontcyclicDelignec}
for related conjecture when $M$ is not assumed to be necessarly formal.
On the contrary, we prove (Corollary~\ref{homologie des lacets egale a cohomologie de Hochschild comme G-algebre}) that Conjecture~\ref{conjecture} is satisfied for $M=S^2$
over $\mathbb{F}_2$.

{\em Acknowledgment:}
We wish to thank Ralph Cohen and Stanford Mathematics department for providing a friendly atmosphere
during my six months of ``delegation CNRS''.
We would like also to thank Yves F\'elix for a discussion simplifying the proof of Theorem~\ref{string topology du cercle}.

\section{The Batalin-Vilkovisky algebra structure on $\mathbb{H}_*(LM)$.}
In this section, we recall the definition of the Batalin-Vilkovisky algebra on
$\mathbb{H}_*(LM;{\Bbbk})$ given by Chas and Sullivan~\cite{Chas-Sullivan:stringtop}
over any commutative ring ${\Bbbk}$
and deduce that this Batalin-Vilkovisky algebra $\mathbb{H}_*(LM;{\Bbbk})$
behaves well with respect to change of rings.

We first recall the definition of the loop product following Cohen and Jones~\cite{CohJon:ahtrostringtopology,Cohen-Hess-Voronov:stringtopacyclhom}.
Let $M$ be a closed oriented smooth manifold of dimension $d$.
The inclusion $e:map(S^1\vee S^1,M)\hookrightarrow LM\times LM$ can be viewed as a codimension $d$
embedding between infinite dimension manifolds~\cite[Proposition 5.3]{Stacey:tdiftopoloops}.
Denote by $\nu$ its normal bundle.
Let $\tau_e:LM\times LM\twoheadrightarrow map(S^1\vee S^1,M)^{\nu}$ its Thom-Pontryagin collapse map.
Recall the umkehr (Gysin) map $e_!$ is the composite of $\tau_e$ and the Thom isomorphism:
$$
H_*(LM\times LM;{\Bbbk})\buildrel{H_*(\tau_e;{\Bbbk})}\over\rightarrow H_*(map(S^1\vee S^1,M)^{\nu};{\Bbbk})
\build\rightarrow_\cong^{\cap u_{\Bbbk}} H_{*-d}(map(S^1\vee S^1,M);{\Bbbk})
$$
The Thom isomorphism is given by taking a relative cap product $\cap$ with a
Thom class for $\nu$, $u_{\Bbbk}\in H^{d}(map(S^1\vee S^1,M)^\nu;{\Bbbk})$.
A Thom class with coefficients in $\mathbb{Z}$, $u_{\mathbb{Z}}$, gives rise a Thom class
$u_{\Bbbk}$ with coefficients in ${\Bbbk}$, under the morphism

$$H^{d}(map(S^1\vee S^1,M);\mathbb{Z})\rightarrow H^{d}(map(S^1\vee S^1,M);{\Bbbk})$$
induced by the ring homomorphism
$\mathbb{Z}\rightarrow{\Bbbk}$~\cite[p. 441-2]{Hatcher:algtop}.
So we have the commutative diagram

\xymatrix{
H_*(LM\times LM;\mathbb{Z})\ar[r]^{e_{!}}\ar[d]
& H_{*-d}(map(S^1\vee S^1,M);\mathbb{Z})\ar[d]\\
H_*(LM\times LM;{\Bbbk})\ar[r]^{e_{!}}
& H_{*-d}(map(S^1\vee S^1,M);{\Bbbk})
} 
Let $\gamma:map(S^1\vee S^1,M)\rightarrow LM$ be the map obtained by composing loops.
The loop product is the composite

\begin{multline*}
H_*(LM;{\Bbbk})\otimes H_*(LM;{\Bbbk})\rightarrow H_*(LM\times LM;{\Bbbk})\\
\buildrel{e_!}\over\rightarrow H_{*-d}(map(S^1\vee S^1,M);{\Bbbk})
\buildrel{H_{*-d}(\gamma;{\Bbbk})}\over\rightarrow H_{*-d}(LM;{\Bbbk})
\end{multline*}
So clearly, we have proved
\begin{lem}\label{loop produit et changement d'anneaux}
The morphism of abelian groups $\mathbb{H}_*(LM;\mathbb{Z})\rightarrow \mathbb{H}_*(LM;{\Bbbk})$
induced by $\mathbb{Z}\rightarrow {\Bbbk}$
is a morphism of graded rings.
\end{lem}
Suppose that the circle $S^1$ acts on a topological space $X$. Then we have an action of the algebra
$H_*(S^1)$ on $H_*(X)$, $$H_*(S^1)\otimes H_*(X)\rightarrow H_*(X).$$
Denote by $[S^1]$ the fundamental class of the circle.
Then we define an operator of degree $1$,
$\BV:H_*(X;{\Bbbk})\rightarrow H_{*+1}(X;{\Bbbk})$ which sends $x$ to 
the image of $[S^1]\otimes x$ under the action.
Since $[S^1]^2=0$, $\BV\circ \BV=0$.
The following lemma is obvious.
\begin{lem}\label{B et changement d'anneaux}
Let $X$ be a $S^1$-space.
We have the commutative diagram

\xymatrix{
H_*(X;\mathbb{Z})\ar[r]^{\BV}\ar[d]
&H_{*+1}(X;\mathbb{Z})\ar[d]\\
H_*(X;{\Bbbk})\ar[r]^{\BV}
&H_{*+1}(X;{\Bbbk})
}
\noindent where the vertical maps are induced by the ring homomorphism $\mathbb{Z}\rightarrow {\Bbbk}$.
\end{lem}
The circle $S^1$ acts on the free loop space on $M$ by rotating the loops.
Therefore we have a operator $\BV$ on $\mathbb{H}_*(LM)$.
Chas and Sullivan~\cite{Chas-Sullivan:stringtop} have showed that
$\mathbb{H}_*(LM)$ equipped with the loop product and the $\BV$ operator,
is a Batalin-Vilkovisky algebra.
\begin{defin}\label{definition algebre de Batalin-Vilkovisky}
A {\it Batalin-Vilkovisky algebra} is a commutative graded algebra
$A$ equipped with an operator $\BV:A\rightarrow A$ of degree $1$
such that $\BV\circ \BV=0$ and
\begin{multline}\label{2nd order derivation}
\BV(abc)=\BV(ab)c+(-1)^{\vert a\vert} a\BV(bc)+(-1)^{(\vert a\vert -1)\vert b\vert}b\BV(ac)\\
-(\BV a)bc-(-1)^{\vert a\vert}a(\BV b)c
-(-1)^{\vert a\vert +\vert b\vert} ab(\BV c).
\end{multline}
\end{defin}
Consider the bracket $\{\;,\;\}$ of degree $+1$ defined
by 
$$\{a,b\}=(-1)^{\vert a\vert}\left(\BV(ab)-(\BV a)b-(-1)^{\vert
  a\vert}a(\BV b)\right)$$
for any $a$, $b\in A$. (\ref{2nd order derivation}) is equivalent to the following relation called the {\it Poisson relation}:
\begin{equation}\label{Poisson relation}
\{a,bc\}=\{a,b\}c+(-1)^{(\vert a\vert+1)\vert b\vert}b\{a,c\}.
\end{equation}
Getzler~\cite[Proposition 1.2]{Getzler:BVAlg} has shown that the $\{\;,\;\}$ is a Lie
  bracket
and therefore that a Batalin-Vilkovisky algebra is a Gerstenhaber algebra.
\begin{defin}\label{definition algebre de Gerstenhaber}

A {\it Gerstenhaber algebra} is a
commutative graded algebra $A$
equipped with a linear map
$\{-,-\}:A \otimes A \to A$ of degree $1$
such that:

\noindent a) the bracket $\{-,-\}$ gives $A$ a structure of graded
Lie algebra of degree $1$. This means that for each $a$, $b$ and $c\in A$

$\{a,b\}=-(-1)^{(\vert a\vert+1)(\vert b\vert+1)}\{b,a\}$ and 

$\{a,\{b,c\}\}=\{\{a,b\},c\}+(-1)^{(\vert a\vert+1)(\vert b\vert+1)}
\{b,\{a,c\}\}.$

\noindent b)  the product and the Lie bracket satisfy the Poisson relation~(\ref{Poisson relation}).
\end{defin}
\noindent Using Lemma~\ref{loop produit et changement d'anneaux} and Lemma~\ref{B et changement d'anneaux}, we deduce
\begin{proposition}\label{BV algebre et changement d'anneaux}
The ${\Bbbk}$-linear map $$\mathbb{H}_*(LM;\mathbb{Z})\otimes_{\mathbb{Z}} {\Bbbk}\hookrightarrow \mathbb{H}_*(LM;{\Bbbk})$$
is an inclusion of Batalin-Vilkovisky algebras.
\end{proposition}
In particular, by the universal coefficient theorem, 
$$\mathbb{H}_*(LM;\mathbb{Z})\otimes_{\mathbb{Z}} {\mathbb{Q}}\cong\mathbb{H}_*(LM;{\mathbb{Q}}).$$
More generally, this Proposition tell us that if $\mbox{Tor}^{\mathbb{Z}}(\mathbb{H}_*(LM;\mathbb{Z}),{\Bbbk})=0$
then the Batalin-Vilkovisky algebra $\mathbb{H}_*(LM;\mathbb{Z})$ determines the Batalin-Vilkovisky algebra
$\mathbb{H}_*(LM;{\Bbbk})$.
\section{The circle and an useful Lemma.}
In this section, we compute the structure of the
Batalin-Vilkovisky algebra on the
homology of the free loop space on the circle $S^1$ using a Lemma
which gives information on the image of $\BV$ on elements of lower degree
in $H_*(LM)$.

\begin{theor}\label{string topology du cercle}
As Batalin-Vilkovisky algebras, the homology of the free loop space on the
circle is given by
$$\mathbb{H}_*(LS^1;\Bbbk)\cong\Bbbk[\mathbb{Z}]\otimes \Lambda a_{-1}.
$$
Denote by $x$ a generator of $\mathbb{Z}$. The operator $\BV$ is
$$\BV(x^i\otimes a_{-1})=i(x^i\otimes 1),\quad\BV(x^i\otimes 1)=0$$
for all $i\in\mathbb{Z}$.
\end{theor}
Let $X$ be a pointed topological space.
Consider the free loop fibration $\Omega X\buildrel{j}\over\hookrightarrow LX
\buildrel{ev}\over\twoheadrightarrow X$.
Denote by $hur_X:\pi_n(X)\rightarrow H_n(X)$ the Hurewicz map.
\begin{lem}\label{Delta, hurewicz et adjonction}
Let $n\in\mathbb{N}$. Let $f\in\pi_{n+1}(X)$.
Denote by $\tilde{f}\in\pi_{n}(\Omega X)$ the adjoint of $f$.
Then $$\left(H_*(ev)\circ \BV\circ H_*(j)\circ hur_{\Omega X}\right)
(\tilde{f})=hur_X(f).$$
\end{lem}
\begin{proof}
Take in homology the image of $[S^1]\otimes [S^n]$ in the following
commutative diagram
$$
\xymatrix{
S^1\times\Omega X\ar[r]^{S^1\times j}
&S^1\times LX\ar[r]^{act_{LX}}
&LX\ar[d]^{ev}\\
S^1\times S^n\ar[r]\ar[u]^{S^1\times\tilde{f}}
&S^1\wedge S^n\ar[r]_f
&X
}
$$
where $act_{LX}:S^1\times LX\rightarrow LX$ is the action of the
circle on $LX$.
\end{proof}
\begin{proof}[Proof of Theorem~\ref{string topology du cercle}]
More generally, let $G$ be a compact Lie group. Consider the homeomorphism
$\Theta_G:\Omega G\times G\buildrel{\cong}\over\rightarrow LG$
which sends the couple $(w,g)$ to the free loop $t\mapsto w(t)g$.
In fact, $\Theta_G$ is an isomorphism of fiberwise monoids.
Therefore by~\cite[part 2) of Theorem 8.2]{gruher-salvatore:genstringtopop},
$$
\mathbb{H}_*(\Theta_G):H_*(\Omega G)\otimes\mathbb{H}_*(G)
\rightarrow\mathbb{H}_*(LG)  
$$
is a morphism of graded algebras. Since $H_*(S^1)$ has no torsion,
$$
\mathbb{H}_*(\Theta_{S^1}):H_*(\Omega S^1)\otimes\mathbb{H}_*(S^1)
\cong \mathbb{H}_*(LS^1)  
$$
is an isomorphism of algebras.
Since $\BV$ preserve path-connected components,
$$
\BV(x^i\otimes a_{-1})=\alpha(x^i\otimes 1)
$$
where $\alpha\in\Bbbk$.
Denote by $\varepsilon_{\Bbbk[\mathbb{Z}]}$ is the canonical augmentation
of the group ring $\Bbbk[\mathbb{Z}]$.
Since $H_*(ev\circ\Theta_{S^1})=
\varepsilon_{\Bbbk[\mathbb{Z}]}\otimes H_*(S^1)$,
$$(H_*(ev)\circ \BV)(x^i\otimes a_{-1})=\alpha 1.$$
On the other hand, applying Lemma~\ref{Delta, hurewicz et adjonction},
to the degree $i$ map $S^1\rightarrow S^1$,
we obtain that $(H_*(ev)\circ \BV\circ H_*(j))(x^i)=i 1$.
Therefore $\alpha=i$.
\end{proof}
\section{Computations using Hochschild homology.}\label{calcul avec l'homologie de Hochschild}
In this section, we compute the Batalin-Vilkovisky algebra
$\mathbb{H}_*(LS^n)$, $n\geq 2$, using the following elementary technique:

The algebra structure has been computed by Cohen, Jones and Yan
using the Serre spectral sequence~\cite{CohJonYan:loophasps}.
On the other hand, the action of $H_*(S^1)$ on $H_*(LS^n)$
can be computed using Hochschild homology.
Using the compatibility between the product and $\BV$,
we determine the Batalin-Vilkovisky algebra $\mathbb{H}_*(LS^n)$
up to isomorphisms. This elementary technique will fail for $\mathbb{H}_*(LS^2)$.

Let $A$ be an augmented differential graded algebra.
Denote by $s\overline{A}$ the suspension of the augmentation ideal
$\overline{A}$, $(s\overline{A})_i=\overline{A}_{i-1}$.
Let $d_1$ be the differential on the tensor product of complexes
$A\otimes T(s\overline{A})$.
The (normalized) Hochschild chain complex, denoted $\mathcal{C}_*(A;A)$, 
is the complex ($A\otimes T(s\overline{A}),d_1+d_2)$
where
\begin{align*}
d_2a[sa_1|\cdots|sa_k]=&(-1)^{\vert a\vert } aa_1[sa_2|\cdots|sa_k]\\
&+\sum _{i=1}^{k-1} (-1)^{\varepsilon_i}{a[sa_1|\cdots|sa_ia_{i+1}|\cdots|sa_k]}\\
&-(-1)^{\vert sa_k\vert\varepsilon_{k-1}} a_ka[sa_1|\cdots|sa_{k-1}];
\end{align*}
Here $\varepsilon_i=\vert a\vert +\vert sa_1\vert+\cdots +\vert sa_i\vert$.

Connes boundary map $B$ is the map of degree $+1$
$$B:A\otimes (s\overline{A})^{\otimes p}\rightarrow
A\otimes (s\overline{A})^{\otimes p+1}$$
defined by

$$B(a_o[sa_1\vert\dots\vert sa_p])=
\sum_{i=0}^p (-1)^{\vert sa_0\dots sa_{i-1}\vert\vert sa_i\dots sa_p\vert}
[sa_i\vert\dots\vert sa_p\vert sa_0\vert\dots\vert sa_{p-1}].
$$

Up to the isomorphism $s^p(A^{\otimes(p+1)})\rightarrow 
A\otimes (sA)^{\otimes p}$, 
$ s^p(a_0[a_1\vert\dots\vert a_p])\mapsto
(-1)^{p\vert a_0\vert+(p-1)\vert a_1\vert+\dots+\vert a_{p-1}\vert}
a_0[sa_1\vert\dots\vert sa_p]$
,
our signs coincides with those of~\cite{VigueM:dechcadgc}.

The Hochschild homology of $A$ (with coefficient in $A$) is the homology
of the Hochschild chain complex:
$$HH_*(A;A):=H_*(\mathcal{C}_*(A;A)).$$
The Hochschild cohomology of $A$ (with coefficient in $A^\vee$)
is the homology of the dual of the Hochschild chain complex:
$$HH^*(A;A^\vee):=H_*(\mathcal{C}_*(A;A)^\vee).$$
Consider the dual of Connes boundary map,
$B^\vee(\varphi)=(-1)^{\vert\varphi\vert}\varphi\circ B$.
On $HH^*(A;A^\vee)$, $B^\vee$ defines an action of $H_*(S^1)$.
\begin{ex}\label{calcul bord de Connes}
Let $n\geq 2$. Let ${\Bbbk}$ be any commutative ring.
Let $A:=H^*(S^n)=\Lambda x_{-n}$ be the exterior algebra on a generator
of lower degree $-n$.
Denote by $[sx]^k:=1[sx\vert\dots\vert sx]$ and
$x[sx]^k:=x[sx\vert\dots\vert sx]$ the elements of $\mathcal{C}_*(A;A)$
where the term $sx$ appears $k$ times. These elements form a basis
of $\mathcal{C}_*(A;A)$. Denote by $[sx]^{k\vee}$, $x[sx]^{k\vee}$, $k\geq 0$,
the dual basis.
The differential $d^\vee$ on $\mathcal{C}_*(A;A)^\vee$ is given by
$d^\vee([sx]^{k\vee})=0$ and
$d^\vee(x[sx]^{k\vee})=\pm\left(1-(-1)^{k(n+1)}\right)[sx]^{(k+1)\vee}$.
The dual of Connes boundary map $B^\vee$ is given by
$$
B^\vee([sx]^{k\vee})=
\begin{cases}
(-1)^{n+1} k\;x[sx]^{(k-1)\vee} & \text{if $(k+1)(n+1)$ is even,}\\
0 & \text{if $(k+1)(n+1)$ is odd}
\end{cases}
$$
and $B^\vee(x[sx]^{k\vee})=0$. We remark that $[sx]^{k\vee}$ is of (lower)
degree $k(n-1)$ and $x[sx]^{k\vee}$ of degree $n+k(n-1)$.
\end{ex}
\begin{theor}~\cite{JonesJ:Cycheh}
Let $X$ be a simply connected space such that $H_*(X;\Bbbk)$
is of finite type in each degree. Then there is a natural isomorphism
of $H_*(S^1)$-modules between the homology of the free loop space
on $X$ and the Hochschild cohomology of the algebra of singular cochain
$S^*(X;\Bbbk)$:
\begin{equation}\label{Jones}
H_*(LX)\cong HH^*(S^*(X;\Bbbk);S^*(X;\Bbbk)^\vee).
\end{equation}
\end{theor}
In this paper, when we will apply this theorem,
$H_*(X;\Bbbk)$ is assumed to be $\Bbbk$-free of finite type in each degree
and $X$ will be always
$\Bbbk$-formal:
the algebra $S^*(X;\Bbbk)$ will be linked by quasi-isomorphisms of
cochain algebras to $H_*(X;\Bbbk)$.
Therefore
\begin{equation}\label{invariance homologie de Hochschild}
HH^*(S^*(X;\Bbbk);S^*(X;\Bbbk)^\vee)
\cong HH^*(H^*(X;\Bbbk);H^*(X;\Bbbk)^\vee).
\end{equation}
\begin{theor}\label{string topology des spheres impairs}
For $n>1$ odd, as Batalin-Vilkovisky algebras,
$$\mathbb{H}_*(LS^n;\Bbbk)=\Bbbk[u_{n-1}]\otimes\Lambda a_{-n},
$$
$$\BV( u_{n-1}^i\otimes a_{-n})=i(u_{n-1}^{i-1}\otimes 1),$$
$$\BV( u_{n-1}^i\otimes 1)=0.$$
\end{theor}
\begin{proof}
As algebras, Cohen, Jones and Yan~\cite{CohJonYan:loophasps} proved that
$\mathbb{H}_*(LS^n;\mathbb{Z})=\Bbbk[u_{n-1}]\otimes\Lambda a_{-n}$
when $\Bbbk=\mathbb{Z}$. Their proof works over any $\Bbbk$
(alternatively, using Proposition~\ref{BV algebre et changement d'anneaux}, we
could assume that $\Bbbk=\mathbb{Z}$).
Computing Connes boundary map on $HH^*(H^*(S^n);H_*(S^n))$
(Example~\ref{calcul bord de Connes}),
we see that $\BV$ on $\mathbb{H}_*(LS^n;\Bbbk)$
is null in even degree and in degree $-n$, and is an isomorphism
in degree $-1$.
Therefore $\BV (u_{n-1}^i\otimes 1)=0$,
$\BV(1\otimes a_{-n})=0$ and
$\BV( u_{n-1}\otimes a_{-n})=\alpha 1$
where $\alpha$ is invertible in $\Bbbk$.
Replacing  $a_{-n}$ by $\frac{1}{\alpha}a_{-n}$ or
$u_{n-1}$ by $\frac{1}{\alpha}u_{n-1}$, we can assume up to isomorphisms
that $\BV( u_{n-1}\otimes a_{-n})=1$.
Therefore $\{u_{n-1},a_{-n}\}=1$. Using the Poisson
relation~(\ref{Poisson relation}), 
$\{u_{n-1}^i,a_{-n}\}=iu_{n-1}^{i-1}$.
Therefore $\BV(u_{n-1}^i\otimes a_{-n})=i(u_{n-1}^{i-1}\otimes 1).$
\end{proof}
\begin{theor}\label{string topology spheres paires sur les entiers}
For $n\geq 2$ even, there exists a constant $\varepsilon_0\in\mathbb{F}_2$ such that
as Batalin-Vilkovisky algebra,
\begin{multline*}
\mathbb{H}_*(LS^n;\mathbb{Z})=\Lambda b\otimes\frac{\mathbb{Z}[a,v]}{(a^2,ab,2av)}\\
=\bigoplus_{k=0}^{+\infty}\mathbb{Z}v_{2(n-1)}^k
\oplus\bigoplus_{k=0}^{+\infty}\mathbb{Z}b_{-1}v^k\oplus\mathbb{Z}a_{-n}
\oplus\bigoplus_{k=1}^{+\infty}\frac{\mathbb{Z}}{2\mathbb{Z}}av^k
\end{multline*}
with $\forall k\geq 0$, $\BV(v^k)=0$, $\BV(av^k)=0$ and
$$\BV(bv^k)=\begin{cases}
(2k+1)v^k+\varepsilon_0 av^{k+1} & \mbox{if }n=2\\
(2k+1)v^k & \mbox{if }n\geq 4.
\end{cases}$$
\end{theor}
\begin{proof}
As algebras, Cohen, Jones and Yan~\cite{CohJonYan:loophasps}
proved the equality.
Computing Connes boundary map on $HH^*(H^*(S^n);H_*(S^n))$
(Example~\ref{calcul bord de Connes}),
we see that $\BV $ on $\mathbb{H}_*(LS^n;\Bbbk)$
is null in even degree and is injective 
in odd degree.

Case $n\neq 2$: this case is simple, since  all the generators
of $\mathbb{H}_*(LS^n)$, $v^k$, $bv^k$ and $av^k$, $k\geq 0$,
have different degrees.
Using Example~\ref{calcul bord de Connes}, we also see that for all $k\geq 0$,
$$
\BV :\mathbb{H}_{-1+2k(n-1)}=\mathbb{Z}b_{-1}v^k
\hookrightarrow \mathbb{H}_{2k(n-1)}=\mathbb{Z}v^k
$$
has cokernel isomorphic to $\frac{\mathbb{Z}}{(2k+1)\mathbb{Z}}$.
Therefore $\BV (bv^k)=\pm(2k+1)v^k$.
By replacing $b_{-1}$ by $-b_{-1}$, we can assume up to isomorphims that
$\BV (b)=1$.
Let $k\geq 1$.
Let $\alpha_k\in\{-2k-1,2k+1\}$ such that $\BV (bv^k)=\alpha_k v^k$.
Using formula~(\ref{2nd order derivation}), we obtain that
$\BV (bv^kv^k)=(2\alpha_k-1) v^{2k}$.
We know that $\BV (bv^{2k})=\pm (4k+1) v^{2k}$.
Therefore $\alpha_k$ must be equal to $2k+1$.

Case $n=2$: this case is complicated, since for $k\geq 0$,
$v^k$ and $av^{k+1}$ have the same degree.
Using Example~\ref{calcul bord de Connes}, we also see that
$$
\BV :\mathbb{H}_{-1+2k}=\mathbb{Z}b_{-1}v^k
\hookrightarrow \mathbb{H}_{2k}=\mathbb{Z}v^k\oplus\frac{\mathbb{Z}}{2\mathbb{Z}}av^{k+1}
$$
has cokernel, denoted $\mbox{Coker}\BV $, isomorphic to $\frac{\mathbb{Z}}{(2k+1)\mathbb{Z}}\oplus\frac{\mathbb{Z}}{2\mathbb{Z}}$.
There exists unique $\alpha_k\in\mathbb{Z}^{*}$ and $\varepsilon_k\in\frac{\mathbb{Z}}{2\mathbb{Z}}$
such that $\BV (bv^k)=\alpha_k v^k+\varepsilon_k av^{k+1}$.
The injective map $\BV $ fits into the commutative diagram of short exact sequences (Noether's Lemma)

\xymatrix{
&0\ar[d]&0\ar[d]&0\ar[d]\\
0\ar[r]
&\mathbb{H}_{-1+2k}\ar[r]^{id}\ar[d]^{\times 2}
&\mathbb{H}_{-1+2k}\ar[r]\ar[d]
&0\ar[r]\ar[d]
&0\\
0\ar[r]
&\mathbb{H}_{-1+2k}\ar[r]^{\BV }\ar[d]
&\mathbb{H}_{2k}\ar[r]\ar[d]
&\mbox{Coker}\BV \ar[r]\ar[d]^{\cong}
&0\\
0\ar[r]
&\frac{\mathbb{Z}}{2\mathbb{Z}}\ar[r]^{\overline{\BV }}\ar[d]
&\frac{\mathbb{Z}}{2\alpha_k\mathbb{Z}}\oplus\frac{\mathbb{Z}}{2\mathbb{Z}}\ar[r]\ar[d]
&\mbox{Coker}\overline{\BV }\ar[r]\ar[d]
&0\\
&0&0&0
}
The cokernel of $\overline{\BV }$, denoted $\mbox{Coker}\overline{\BV }$ is of cardinal $2\vert\alpha_k\vert$.
So $\vert\alpha_k\vert=2k+1$. Therefore
$\BV (bv^k)=\pm (2k+1) v^k+\varepsilon_k av^{k+1}$.

By replacing $b_{-1}$ by $-b_{-1}$, we can assume up to isomorphims that
$\BV (b)=1+\varepsilon_0 av$.
Using formula~(\ref{2nd order derivation}), we obtain that
$$\BV (bv^k v^l)=(\alpha_k+\alpha_l-1) v^{k+l}+(\varepsilon_k+\varepsilon_l-\varepsilon_0)av^{k+l+1}.$$
Therefore
$$
\BV (bv^kv^k)=(2\alpha_k-1)v^{2k}+\varepsilon_0 av^{2k+1}=\pm(4k+1)v^{2k}+\varepsilon_{2k}av^{2k+1}.
$$
So $\alpha_k=2k+1$, $\varepsilon_{2k}=\varepsilon_{0}$ and
$\varepsilon_{2k+1}=\varepsilon_{2k}+\varepsilon_1-\varepsilon_0=\varepsilon_1$.

The map $\Theta:\mathbb{H}_*(LS^2)\rightarrow\mathbb{H}_*(LS^2)$
given by
$\Theta(b_{-1}v^k)=b_{-1}v^k$,
$\Theta(v^k)=v^k+kav^{k+1}$,
$\Theta(av^k)=av^k$, $k\geq 0$
is an involutive isomorphism of algebras.
Therefore, by replacing $v$ by $v+av^2$, we can assume that $\varepsilon_1=\varepsilon_0$.
So we have proved
$$\BV (bv^k)=(2k+1)v^k+\varepsilon_0av^{k+1},\quad k\geq 0.$$\end{proof}
These two cases $\varepsilon_0=0$ and $\varepsilon_0=1$
correspond to two non-isomorphic Batalin-Vilkovisky algebras whose underlying Gerstenhaber algebras are the same.
Therefore even if we have not yet computed the Batalin-Vilkovisky algebra
$\mathbb{H}_*(LS^2;\mathbb{Z})$, we have computed its underlying Gerstenhaber algebra. Using the definition of the
bracket, straightforward computations give the following corollary.
\begin{cor}
For $n\geq 2$ even, as Gerstenhaber algebra
$$\mathbb{H}_*(LS^n;\mathbb{Z})=\Lambda b_{-1}\otimes\frac{\mathbb{Z}[a_{-n},v_{2(n-1)}]}{(a^2,ab,2av)}$$
with $\{v^k,v^l\}=0$, $\{bv^k,v^l\}=-2lv^{k+l}$, $\{bv^k,bv^l\}=2(k-l)bv^{k+l}$,
$\{a,v^l\}=0$, $\{av^k,bv^l\}=-(2l+1)av^{k+l}$ and $\{av^k,av^l\}=0$ for all $k,l\geq 0$.
\end{cor}
\section{When Hochschild cohomology is a Batalin-Vilkovisky algebra}
In this section, we recall the structure of Gerstenhaber algebra
on the Hochschild cohomology of an algebra whose degrees are bounded.
We recall from~\cite{Tradler:bvalgcohiiip,MenichiL:BValgaccoHa,Tradler-Zeinalian:ontcyclicDelignec,Kaufmann:procvDelignecvc}
the Batalin-Vilkovisky algebra on the Hochschild cohomology of the cohomology
$H^*(M)$ of a closed oriented manifold $M$.
We compute this Batalin-Vilkovisky algebra
$HH^*(H^*(M);H^*(M))$ when $M$ is a sphere.

Through this section, we will work over the prime field $\mathbb{F}_2$.
Let $A$ be an augmented graded algebra such that the augmentation ideal
$\overline{A}$ is concentrated in degree $\leq -2$ and bounded below
(or concentrated in degree $\geq 0$ and bounded above).
Then the (normalized) Hochschild cochain complex, denoted $\mathcal{C}^{*}(A,A)$, is the complex
$$\mbox{Hom}(Ts\overline{A},A)\cong \oplus_{p\geq 0}\mbox{Hom}((s\overline{A})^{\otimes p},A)$$
with a differential $d_2$.
For $f\in\mbox{Hom}((s\overline{A})^{\otimes p},A)$,
the differential $d_2f\in\mbox{Hom}((s\overline{A})^{\otimes p+1},A)$ is
given by
\begin{multline*}
(d_2f)([sa_1|\cdots|sa_{p+1}]):=
a_1f([sa_2|\cdots|sa_{p+1}])\\
+\sum_{i=1}^pf([sa_1|\cdots|s(a_ia_{i+1})|\cdots|sa_{p+1}])
+f([sa_1|\cdots|sa_{p}])a_p
\end{multline*}
The Hochschild cohomology of $A$ with coefficient in $A$ is the homology
of the Hochschild cochain complex:
$$HH^*(A;A):=H_*(\mathcal{C}^*(A;A)).$$
We remark that $HH^*(A;A)$ is bigraded. Our degree is sometimes called the total degree: sum of the external degree
and the internal degree.
The Hochschild cochain complex $\mathcal{C}^{*}(A,A)$ is a differential graded algebra.
For $f\in\mbox{Hom}((s\overline{A})^{\otimes p},A)$ and  $g\in\mbox{Hom}((s\overline{A})^{\otimes q},A)$,
the (cup) product of $f$ and $g$, $f\cup g\in\mbox{Hom}((s\overline{A})^{\otimes p+q},A)$ is defined by
$$(f\cup g)([sa_1|\cdots|sa_{p+q}]):=f([sa_1|\cdots|sa_{p}])g([sa_{p+1}|\cdots|sa_{p+q}]).$$
The Hochschild cochain complex $\mathcal{C}^{*}(A,A)$ has also a Lie bracket of (lower) degree $+1$.
\begin{multline*}
(f\overline{\circ} g)([sa_1|\cdots|sa_{p+q-1}]):=\\
\sum_{i=1}^p f\left([sa_1|\cdots|sa_{i-1}|sg([sa_{i}|\cdots|sa_{i+q-1}])|sa_{i+q}|\cdots|sa_{p+q-1}]\right).$$
\end{multline*}
$\{f,g\}=f\overline{\circ} g-g\overline{\circ} f$.
Our formulas are the same as in the non graded case~\cite{Gerstenhaber:cohosar}.
We remark that if $A$ is not assumed to be bounded,
the formulas are more complicated.
Gerstenhaber has showed that $HH^*(A;A)$ equipped with the cup product and the Lie bracket is a Gerstenhaber algebra.

Let $M$ be a closed $d$-dimensional smooth manifold.
Poincar\'e duality induces an isomorphism of $H^*(M;\mathbb{F}_2)$-modules of (lower) degree $d$.

\begin{equation}\label{dualite de Poincare}
\Theta:H^*(M;\mathbb{F}_2)\buildrel{\cap [M]}\over\rightarrow H_*(M;\mathbb{F}_2)\cong H^*(M;\mathbb{F}_2)^\vee.
\end{equation}
More generally, let $A$ be a graded algebra equipped with an isomorphism of $A$-bimodules of degree $d$,
$\Theta:A\buildrel{\cong}\over\rightarrow A^{\vee}$.
Then we have the isomorphism
$$
HH^*(A,\Theta):HH^{*}(A,A)\buildrel{\cong}\over\rightarrow HH^*(A,A^{\vee}).
$$
Therefore on $HH^{*}(A,A)$, we have both a Gerstenhaber algebra structure and an operator $\BV$ given by the
dual of Connes boundary map $B$.
Motivated by the Batalin-Vilkovisky algebra structure of Chas-Sullivan on $\mathbb{H}_*(LM)$, Thomas Tradler~\cite{Tradler:bvalgcohiiip} proved that
$HH^{*}(A,A)$ is a Batalin-Vilkovisky algebra.
See~\cite[Theorem 1.6]{MenichiL:BValgaccoHa} for an explicit proof.
In \cite{Kaufmann:procvDelignecvc} or
\cite[Corollary 3.4]{Tradler-Zeinalian:ontcyclicDelignec} or
\cite[Section 1.4]{Costello:tcftcalabiyaucat} or
\cite[Theorem B]{Kaufmann:modsahcfaco} or
\cite[Section 11.6]{Kontsevich-Soibelman:noainfinityalg},
this Batalin-Vilkovisky algebra structure on $HH^{*}(A,A)$ extends
to a structure of algebra on the Hochschild cochain complex
$\mathcal{C}^{*}(A,A)$ over various operads or PROPs: the so-called
cyclic Deligne conjecture.
Let us compute this Batalin-Vilkovisky algebra structure when $M$ is a sphere.
\begin{proposition}\label{Cohomologie de Hochschild des spheres modulo 2}
(\cite{Westerland:DyerLoitstosaps}
and~\cite[Corollary 4.2]{Westerland:stringhomologyspheres})
Let $d\geq 2$. As Batalin-Vilkovisky algebra, the Hochschild cohomology of $H^*(S^d;\mathbb{F}_2)=\Lambda x_{-d}$,
$$HH^*(H^*(S^d;\mathbb{F}_2);H^*(S^d;\mathbb{F}_2))\cong\Lambda g_{-d}\otimes \mathbb{F}_2[f_{d-1}]$$
with $\BV(g_{-d}\otimes f_{d-1}^k)=k(1\otimes f_{d-1}^{k-1})$ and $\BV(1\otimes f_{d-1}^k)=0, k\geq 0$.
In particular, the underlying Gerstenhaber algebra is given by
$\{f^k,f^l\}=0$, $\{gf^k,f^l\}=lf^{k+l-1}$ and $\{gf^k,gf^l\}=(k-l)gf^{k+l-1}$ for $k$, $l\geq 0$.
\end{proposition}
\begin{proof}
Denote by $A:=H^*(S^d;\mathbb{F}_2)$.
The differential on $\mathcal{C}^*(A;A)$ is null.
Let $f\in\mbox{Hom}(s\overline{A},A)\subset \mathcal{C}^*(A;A)$ such that $f([sx])=1$.
Let $g\in\mbox{Hom}(\mathbb{F}_2,A)=\mbox{Hom}((s\overline{A})^{\otimes 0},A)\subset \mathcal{C}^*(A;A)$ such that $g([])=x$.
The $k$-th power of $f$ is the map $f^k\in\mbox{Hom}((s\overline{A})^{\otimes k},A)$ such that
$f^k([sx|\cdots|sx])=1$.
The cup product $g\cup f^k\in\mbox{Hom}((s\overline{A})^{\otimes k},A)$ 
sends $[sx|\cdots|sx]$ to $x$.
So we have proved that $\mathcal{C}^*(A;A)$ is isomorphic to the tensor product of graded
algebras $\Lambda g_{-d}\otimes \mathbb{F}_2[f_{d-1}]$.

The unit $1$ and $x_{-d}$ form a linear basis of $H^*(S^d)$.
Denote by $1^{\vee}$ and $x^{\vee}$ the dual basis of $A^\vee=H^*(S^d)^\vee$.
Poincar\'e duality induces the isomorphism $\Theta:H^*(S^d)\buildrel{\cong}\over\rightarrow H^*(S^d)^\vee$, $1\mapsto x^\vee$ and
$x\mapsto 1^\vee$.
The two families of elements of the form $1[sx|\cdots|sx]$ and of the form $x[sx|\cdots|sx]$ forms a basis
of $\mathcal{C}_*(A;A)$. Denote by $1[sx|\cdots|sx]^\vee$ and $x[sx|\cdots|sx]^\vee$ the dual basis in $\mathcal{C}_*(A;A)^\vee$.
The isomorphism $\Theta$ induces an isomorphism of complexes of degree $d$,
$\widehat\Theta:\mathcal{C}^*(A;A)\build\rightarrow_\cong^{\mathcal{C}^*(A;\Theta)}
\mathcal{C}^*(A;A^\vee)\buildrel{\cong}\over\rightarrow\mathcal{C}_*(A;A)^\vee$.
Explicitly~\cite[Section 4]{MenichiL:BValgaccoHa} this isomorphism sends
$f\in\mbox{Hom}((s\overline{A})^{\otimes p},A)$ to the linear map
$\widehat\Theta(f)\in (A\otimes (s\overline{A})^{\otimes p})^\vee\subset\mathcal{C}_*(A;A)^\vee$ defined
by
$$\widehat\Theta(f)(a_0[sa_1|\cdots|sa_p])=\left((\Theta\circ f)[sa_1|\cdots|sa_p]\right)(a_0).$$
Here with $A=\Lambda x$,
$\widehat\Theta(f^k)=x[sx|\cdots|sx]^\vee$
and $\widehat\Theta(g\cup f^k)=1[sx|\cdots|sx]^\vee$.
Computing Connes boundary map $B^\vee$ on $\mathcal{C}_*(A;A)^\vee$
(Example~\ref{calcul bord de Connes}) and using that by definition of $\BV$,
$\widehat\Theta\circ\BV=B^\vee\circ\widehat\Theta$, we obtain the desired formula for $\BV$.
\end{proof}
\section{The Gerstenhaber algebra $\mathbb{H}_*(LS^2;\mathbb{F}_2$)}\label{algebre de Gerstenhaber sur la sphere modulo 2}
Using the same Hochschild homology technique as in
section~\ref{calcul avec l'homologie de Hochschild},
we compute up to an indeterminacy,
the Batalin-Vilkovisky algebra $\mathbb{H}_*(LS^2;\mathbb{F}_2)$.
Nevertheless, this will give the complete description of
the underlying Gerstenhaber algebra on $\mathbb{H}_*(LS^2;\mathbb{F}_2)$.
\begin{lem}\label{BV algebre des lacets sur la sphere a epsilon pres modulo 2}
There exist a constant $\varepsilon\in\{0,1\}$ such that
as Batalin-Vilkovisky algebra, the homology of the free loop on the sphere $S^2$ is
$$\mathbb{H}_*(LS^2;\mathbb{F}_2)=\Lambda a_{-2}\otimes \mathbb{F}_2[u_1],$$
$$\BV(a_{-2}\otimes u_1^k)=k(1\otimes u_1^{k-1}+\varepsilon a_{-2}\otimes u_1^{k+1} )\mbox{ and }\BV(1\otimes u_1^k)=0, k\geq 0.$$
\end{lem}
\begin{proof}
In~\cite{CohJonYan:loophasps}, Cohen, Jones and Yan proved that 
the Serre spectral sequence for the free loop fibration
$\Omega M\buildrel{j}\over\hookrightarrow LM\buildrel{ev}\over\twoheadrightarrow M$
is a spectral sequence of algebras converging toward the algebra $\mathbb{H}_*(LM)$.
Using Hochschild homology, we see that there is an isomorphism of vector spaces
$\mathbb{H}_*(LS^2;\mathbb{F}_2)\cong \mathbb{H}_*(S^2;\mathbb{F}_2)\otimes H_*(\Omega S^2;\mathbb{F}_2)$.
Therefore the Serre spectral sequence collapses.
Since there is no extension problem, we have the isomorphism of algebras
$$\mathbb{H}_*(LS^2;\mathbb{F}_2)\cong \mathbb{H}_*(S^2;\mathbb{F}_2)\otimes H_*(\Omega S^2;\mathbb{F}_2)=
\Lambda (a_{-2})\otimes\mathbb{F}_2[u_1].$$
Computing Connes boundary map on $HH^*(H^*(S^2;\mathbb{F}_2);H_*(S^2;\mathbb{F}_2))$ (Example~\ref{calcul bord de Connes}),
we see that $\BV$ on $\mathbb{H}_*(LS^2;\mathbb{F}_2)$
is null in even degree and that
$$\BV:\mathbb{H}_{2k-1}\rightarrow \mathbb{H}_{2k}$$
is a linear map of rank $1$, $k\geq 0$. In particular $\BV$ is injective in degree~$-1$.

Applying Lemma~\ref{Delta, hurewicz et adjonction}, to the identity map $id:S^2\rightarrow S^2$,
we see that the composite
$$H_1(\Omega S^2;\mathbb{F}_2)\buildrel{H_1(j;\mathbb{F}_2)}\over\rightarrow H_1(LS^2;\mathbb{F}_2)
 \buildrel{\BV}\over\rightarrow  H_2(LS^2;\mathbb{F}_2)
\buildrel{H_2(ev;\mathbb{F}_2)}\over\rightarrow H_2(S^2;\mathbb{F}_2)
$$
is non zero.
Since $\mathbb{H}_*(ev)$ is a morphism of algebras, $\mathbb{H}_{0}(ev)(a_{-2}u_1^2)=0$.
And so $\BV(a_{-2}u_1)=1+\varepsilon a_{-2}u_1^2$ with $\varepsilon\in\mathbb{F}_2$.

We remark that when $b=c$, formula~(\ref{2nd order derivation})
takes the simple form
\begin{equation}\label{B derivation pour produit avec un carre}
\BV(ab^2)=\BV(a)b^2+a\BV(b^2).
\end{equation}
Using this formula, we obtain that
$$\BV(a_{-2}u_1^{2k+1})=\BV((a_{-2}u_1)(u_1^{k})^2)=u_1^{2k}+\varepsilon a_{-2}u_1^{2k+2}\quad k\geq 0.$$
%
Since $\BV:\mathbb{H}_1=\mathbb{F}_2 a_{-2} u_1^3\oplus
\mathbb{F}_2u_1\rightarrow\mathbb{H}_2$ is of rank $1$ and
$\BV(a_{-2}u_1^{3})\neq 0$,
$\BV(u_1)=\lambda \BV(a_{-2}u_1^3)$ with $\lambda=0$ or $\lambda=1$.
Using again formula~(\ref{B derivation pour produit avec un carre}),
we have that
$$\BV(u_1^{2k+1})=\BV(u_1(u_1^{k})^{2})=\lambda\BV(a_{-2}u_1^{3})u_1^{2k}=\lambda\BV(a_{-2}u_1^{2k+3}),k\geq 0.$$
So finally
$$\BV(a_{-2}u_1^k)=ku_1^{k-1}+\varepsilon ka_{-2}u_1^{k+1} \mbox{ and }
\BV(u_1^k)=\lambda\BV(a_{-2}u_1^{k+2}), k\geq 0.$$
The cases $\lambda=0$ and $\lambda=1$ correspond to isomorphic Batalin-Vilkovisky algebras:
Let $\Theta:\mathbb{H}_*(LS^2;\mathbb{F}_2)\rightarrow\mathbb{H}_*(LS^2;\mathbb{F}_2)$
be an automorphism of algebras which is not the identity.
Since $\Theta(a_{-2})\neq 0$, $\Theta(a_{-2})=a_{-2}$.
Since $\Theta(a_{-2})$ and $\Theta(u_1)$ must generate the algebra
$\Lambda a_{-2}\otimes\mathbb{F}_2[u_1]$, $\Theta(u_1)\neq a_{-2}u_1^3$.
Since $\Theta(u_1)\neq u_1$, $\Theta(u_1)=u_1+a_{-2}u_1^3$.
Therefore there is an unique automorphism of algebras
$\Theta:\mathbb{H}_*(LS^2;\mathbb{F}_2)\rightarrow\mathbb{H}_*(LS^2;\mathbb{F}_2)$ which is not the identity.
Explicitly, $\Theta$ is given by
$\Theta(u_1^k)=u_1^k+ka_{-2}u_1^{k+2}$,
$\Theta(a_{-2}u_1^k)=a_{-2}u_1^k$, $k\geq 0$.
One can check that $\Theta$ is an involutive isomorphism of Batalin-Vilkovisky algebras
who transforms the cases $\lambda=0$ into the cases $\lambda=1$ without changing $\varepsilon$.
Therefore, by replacing $u_1$ by $u_1+a_{-2}u_1^3$, we can assume that $\lambda=0$.
\end{proof}
Consider the four Batalin-Vilkovisky algebras $\Lambda a_{-2}\otimes\mathbb{F}_2[u_1]$ with
$\BV(a_{-2}\otimes u_1^k)=k(1\otimes u_1^{k-1}+\varepsilon a_{-2}\otimes u_1^{k+1} )$,
$\BV(1\otimes u_1^k)=\lambda\BV(a_{-2}u_1^{k+2})$, $k\geq 0$, given the different values of
$\varepsilon$, $\lambda\in\{0,1\}$.
These four Batalin-Vilkovisky algebras have only two underlying Gerstenhaber algebras given by
$\{u_1^k,u_1^l\}=0$, $\{a_{-2}u_1^k,u_1^l\}=lu^{k+l-1}+l(\varepsilon-\lambda)a_{-2}u^{k+l+1}$
and $\{a_{-2}u_1^k,a_{-2}u_1^l\}=(k-l)a_{-2}u^{k+l-1}$ for $k$, $l\geq 0$.
Via the above isomorphism $\Theta$, these two Gerstenhaber algebras are isomorphic.
\begin{cor}\label{homologie des lacets egale a cohomologie de Hochschild comme G-algebre}
The free loop space modulo $2$ homology $\mathbb{H}_*(LS^2;\mathbb{F}_2)$ is isomorphic
as Gerstenhaber algebra to the Hochschild cohomology of $H^*(S^2;\mathbb{F}_2)$,
$HH^{*}(H^*(S^2;\mathbb{F}_2);H^*(S^2;\mathbb{F}_2))$.
\end{cor}
\section{The Batalin-Vilkovisky algebra $\mathbb{H}_*(LS^2)$}
In this section, we complete the calculations of the Batalin-Vilkovisky
algebras $\mathbb{H}_*(LS^2;\mathbb{F}_2)$ and $\mathbb{H}_*(LS^2;\mathbb{Z})$
started respectively
in sections~\ref{algebre de Gerstenhaber sur la sphere modulo 2}
and~\ref{calcul avec l'homologie de Hochschild}, using a purely homotopic method.
\begin{theor}\label{homologie des lacets sur la sphere modulo 2}
As Batalin-Vilkovisky algebra, the homology of the free loop space on the sphere $S^2$ with mod $2$ coefficients is
$$\mathbb{H}_*(LS^2;\mathbb{F}_2)=\Lambda a_{-2}\otimes \mathbb{F}_2[u_1],$$
$$\BV(a_{-2}\otimes u_1^k)=k(1\otimes u_1^{k-1}+a_{-2}\otimes u_1^{k+1} )
\mbox{ and }\BV(1\otimes u_1^k)=0, k\geq 0.$$
\end{theor}
\begin{theor}\label{homologie des lacets sur la sphere modulo Z}
With integer coefficients, as Batalin-Vilkovisky algebra,
\begin{multline*}
\mathbb{H}_*(LS^2;\mathbb{Z})=\Lambda b\otimes\frac{\mathbb{Z}[a,v]}{(a^2,ab,2av)}\\
=\bigoplus_{k=0}^{+\infty}\mathbb{Z}v_{2}^k
\oplus\bigoplus_{k=0}^{+\infty}\mathbb{Z}b_{-1}v^k\oplus\mathbb{Z}a_{-2}
\oplus\bigoplus_{k=1}^{+\infty}\frac{\mathbb{Z}}{2\mathbb{Z}}av^k
\end{multline*}
with $\forall k\geq 0$, $\BV(v^k)=0$, $\BV(av^k)=0$ and
$\BV(bv^k)=(2k+1)v^k+av^{k+1}$.
\end{theor}
Denote by $\sectiontriviale:X\hookrightarrow LX$ the trivial section of the
evaluation map $ev:LX\twoheadrightarrow X$.
\begin{lem}\label{images de B de la section modulo 2}
The image of $\BV :H_1(LS^2;\mathbb{F}_2)\rightarrow H_2(LS^2;\mathbb{F}_2)$ is not contained
in the image of $H_2(\sectiontriviale;\mathbb{F}_2):H_2(S^2;\mathbb{F}_2)\hookrightarrow H_2(LS^2;\mathbb{F}_2)$.
\end{lem}
\begin{lem}\label{images de B de la section modulo Z}
The image of $\BV :H_1(LS^2;\mathbb{Z})\rightarrow H_2(LS^2;\mathbb{Z})$ is not contained
in the image of $H_2(\sectiontriviale;\mathbb{Z}):H_2(S^2;\mathbb{Z})\hookrightarrow H_2(LS^2;\mathbb{Z})$.
\end{lem}
\begin{proof}[Proof of Lemma~\ref{images de B de la section modulo Z} assuming Lemma~\ref{images de B de la section modulo 2}]
Consider the commutative diagram

\xymatrix{
H_1(LS^2;\mathbb{Z})\otimes_{\mathbb{Z}}\mathbb{F}_2\ar[r]^{\cong}\ar[d]_{\BV\otimes_{\mathbb{Z}}\mathbb{F}_2}
& H_1(LS^2;\mathbb{F}_2)\ar[d]^{\BV}\\
H_2(LS^2;\mathbb{Z})\otimes_{\mathbb{Z}}\mathbb{F}_2\ar[r]^{\cong}
& H_2(LS^2;\mathbb{F}_2)\\
H_2(S^2;\mathbb{Z})\otimes_{\mathbb{Z}}\mathbb{F}_2\ar[r]^{\cong}
\ar[u]^{H_2(\sectiontriviale;\mathbb{Z})\otimes_{\mathbb{Z}}\mathbb{F}_2}
& H_2(S^2;\mathbb{F}_2)\ar[u]_{H_2(\sectiontriviale;\mathbb{F}_2)}
}
Since $H_1(LS^2;\mathbb{Z})\cong H_0(LS^2;\mathbb{Z})\cong\mathbb{Z}$, the
horizontal arrows are isomorphisms by the universal coefficient theorem.
The top rectangle commutes according Lemma~\ref{B et changement d'anneaux}.

Suppose that the image of $\BV :H_1(LS^2;\mathbb{Z})\rightarrow H_2(LS^2;\mathbb{Z})$ is included
in the image of $H_2(\sectiontriviale;\mathbb{Z})$.
Then the image of $\BV\otimes_{\mathbb{Z}}\mathbb{F}_2$ is included
in the image of $H_2(\sectiontriviale;\mathbb{Z})\otimes_{\mathbb{Z}}\mathbb{F}_2$.
Using the above diagram, the image of $\BV :H_1(LS^2;\mathbb{F}_2)\rightarrow H_2(LS^2;\mathbb{F}_2)$ is included
in the image of $H_2(\sectiontriviale;\mathbb{F}_2)$.
This contradicts Lemma~\ref{images de B de la section modulo 2}.
\end{proof}
\begin{proof}[Proof of Theorem~\ref{homologie des lacets sur la sphere modulo 2} assuming
Lemma~\ref{images de B de la section modulo 2}]
It suffices to show that the constant $\varepsilon$ in Lemma~\ref{BV algebre des lacets sur la sphere a epsilon pres modulo 2}
is not zero.
Suppose that $\varepsilon=0$.
Then by Lemma~\ref{BV algebre des lacets sur la sphere a epsilon pres modulo 2},
$\BV(a_{-2}\otimes u_1)=1$.

It is well known that $\mathbb{H}_*(\sectiontriviale):\mathbb{H}_*(M)\rightarrow\mathbb{H}_*(LM)$ is a morphism of algebras.
In particular, let $[S^2]$ be the fundamental class of $S^2$, $H_2(\sectiontriviale)([S^2])$ is the unit
of $\mathbb{H}_*(LS^2)$.
So $\BV(a_{-2}\otimes u_1)=H_2(\sectiontriviale)([S^2])$. This contradicts
Lemma~\ref{images de B de la section modulo 2}.
\end{proof}
The proof of Theorem~\ref{homologie des lacets sur la sphere modulo Z} assuming
Lemma~\ref{images de B de la section modulo Z} is the same.
To complete the computation of this Batalin-Vilkovisky algebra on the homology of the free loop space
of a manifold, we will relate it to another structure of Batalin-Vilkovisky algebra that arises in
algebraic topology: the homology of the double loop space.

Let $X$ be a pointed topological space.
The circle $S^1$ acts on the sphere $S^2$ by ``rotating the earth''.
Therefore the circle also acts on $\Omega^2 X=map\left((S^2,\mbox{North pole}),(X,*)\right)$.
So we have a induced operator $\BV:H_*(\Omega^2 X)\rightarrow H_{*+1}(\Omega^2 X)$.
With Theorem~\ref{Gaudens Menichi} and the following Proposition, we will able to prove
Lemma~\ref{images de B de la section modulo 2}.
\begin{proposition}\label{retract commute avec les actions du cercle}
Let $X$ be a pointed topological space.
There is a natural morphism 
$r:L\Omega X\rightarrow map_*(S^2,X)$
of $S^1$-spaces between the free loop on the pointed loop of $X$
and the double pointed loop space of $X$ such that:

$\bullet$ If we identify $S^2$ and $S^1\wedge S^1$, $r$ is a retract up to homotopy of the inclusion $j:\Omega(\Omega X)\hookrightarrow L(\Omega X)$,

$\bullet$ The composite
$r\circ \sectiontriviale:\Omega X\hookrightarrow L(\Omega X)\rightarrow map_*(S^2,X)$
is homotopically trivial.
\end{proposition}
\begin{proof}
Let $\sigma:S^2\twoheadrightarrow \frac{S^1\times S^1}{S^1\times *}=S^1_+\wedge S^1$
be the quotient map that identifies the North pole and the South pole on the earth $S^2$.
The circle $S^1$ acts without moving the based point on $S^1_+\wedge S^1$ by multiplication on the first factor.
On the torus $S^1\times S^1$, the circle can act by
multiplication on both factors. But when you pinch a circle to a point in the torus,
the circle can act only on one factor.
If we make a picture, we easily see that
$\sigma:S^2\twoheadrightarrow S^1_+\wedge S^1$
is compatible with the actions of $S^1$.
Therefore $r:=map_*(\sigma,X):L\Omega X\rightarrow map_*(S^2,X)$ is  a morphism
of $S^1$-spaces.

$\bullet$
Let $\pi:S^1_+\wedge S^1\twoheadrightarrow S^1\wedge S^1=\frac{S^1_+\wedge S^1}{*\times S^1}$
be the quotient map. The inclusion map
$j:\Omega(\Omega X)\rightarrow L(\Omega X)$ is $map_*(\pi,X)$.
The composite $\pi\circ\sigma:S^2\twoheadrightarrow S^1\wedge S^1$ is the quotient map
obtained by identify a meridian to a point in the sphere $S^2$.
The composite $\pi\circ\sigma$ can also be viewed as the quotient map from the non reduced
suspension of $S^1$ to the reduced suspension of $S^1$.
So the composite $\pi\circ\sigma:S^2\twoheadrightarrow S^1\wedge S^1$ is a homotopy equivalence.
Let $\Theta:S^1\wedge S^1\buildrel{\cong}\over\rightarrow S^2$ be any given homeomorphism.
The composite $\Theta\circ\pi\circ\sigma:S^2\rightarrow S^2$ is of degree $\pm 1$.
The reflection through the equatorial plane is a morphism of $S^1$-spaces.
By replacing eventually $\sigma$ by its composite with the previous reflection,
we can suppose that $\Theta\circ\pi\circ\sigma:S^2\rightarrow S^2$ is homotopic to the
identity map of $S^2$, i. e. $\sigma\circ\Theta$ is a section of $\pi$ up to homotopy.
Therefore $map_*(\sigma\circ\Theta,X)=map_*(\Theta,X)\circ r$ is a retract of $j$ up to homotopy.

$\bullet$ Let $\rho:S^1_+\wedge S^1=\frac{S^1\times S^1}{S^1\times *}\twoheadrightarrow S^1$
be the map induced by the projection on the second factor.
Since $\pi_2(S^1)=*$, the composite $\rho\circ\sigma$ is homotopically trivial.
Therefore $r\circ\sectiontriviale$, the composite of $r=map_*(\sigma,X)$
and $\sectiontriviale=map_*(\rho, X):\Omega X\rightarrow L(\Omega X)$ is also homotopically trivial. 
\end{proof}
\begin{proof}[Proof of Lemma~\ref{images de B de la section modulo 2}]
Denote by $ad_{S^n}:S^n\rightarrow\Omega S^{n+1}$ the adjoint of the identity map $id:S^{n+1}\rightarrow S^{n+1}$.
The map $L(ad_{S^2}):LS^2\rightarrow L\Omega S^3$ is obviously a morphism of $S^1$-spaces.
Therefore using Proposition~\ref{retract commute avec les actions du cercle},
the composite $r\circ L(ad_{S^2}):LS^2\rightarrow L\Omega S^3\rightarrow \Omega^2 S^3$
is also a morphism of $S^1$-spaces. Therefore $H_*(r\circ L(ad_{S^2}))$ commutes with the corresponding operators $\BV $
in $H_*(LS^2)$ and $H_*(\Omega^2 S^3)$.

Consider the commutative diagram up to homotopy
\begin{equation}~\label{retract et section et ad}
\xymatrix{
\Omega S^2\ar[r]^{j}\ar[d]_{\Omega (ad_{S^2})}
& LS^2\ar[d]_{L(ad_{S^2})}
& S^2\ar[l]_{\sectiontriviale}\ar[d]^{ad_{S^2}}\\
\Omega^2 S^3\ar[r]^{j}\ar[dr]_{id}
& L\Omega S^3\ar[d]^{r}
& \Omega S^3\ar[l]_{\sectiontriviale}\ar[dl]^{*}\\
& \Omega^2 S^3
}
\end{equation}
Using the left part of this diagram, we see that $\pi_1(r\circ L(ad))$ maps the generator of
$\pi_1(LS^2)=\mathbb{Z}(j\circ ad_{S^1})$ to the composite
$\Omega (ad_{S^2})\circ ad_{S^1}:S^1\rightarrow\Omega S^2\rightarrow \Omega^2 S^3$ which is the generator of
$\pi_1(\Omega^2 S^3)\cong\mathbb{Z}$. Therefore $\pi_1(r\circ L(ad))$ is an isomorphism.

So we have the commutative diagram

\xymatrix{
\pi_1(LS^2)\otimes\mathbb{F}_2 \ar[r]^{hur}_{\cong}\ar[d]_{\pi_1(r\circ L(ad_{S^2}))\otimes\mathbb{F}_2}^{\cong}
& H_1(LS^2;\mathbb{F}_2)\ar[r]^{\BV }\ar[d]_{H_1(r\circ L(ad_{S^2});\mathbb{F}_2)}
& H_2(LS^2;\mathbb{F}_2)\ar[d]^{H_2(r\circ L(ad_{S^2});\mathbb{F}_2)}\\
\pi_1(\Omega^2 S^3)\otimes\mathbb{F}_2 \ar[r]^{hur}_{\cong}
& H_1(\Omega^2 S^3;\mathbb{F}_2)\ar[r]^{\BV }
& H_2(\Omega^2 S^3;\mathbb{F}_2)
}
\noindent By Theorem~\ref{Gaudens Menichi}, $\BV : H_1(\Omega^2 S^3;\mathbb{F}_2)\rightarrow
H_2(\Omega^2 S^3;\mathbb{F}_2)$ is non zero.
Therefore using the above diagram, the composite $H_2(r\circ L(ad_{S^2}))\circ \BV $
is also non zero.
On the other hand, using the right part of diagram~(\ref{retract et
  section et ad}), we have that the composite
$H_2(r\circ L(ad_{S^2}))\circ H_2(\sectiontriviale)$
is null.
\end{proof}
\begin{cor}\label{homologie des lacets different de cohomologie de Hochschild comme BV-algebre}
The free loop space modulo $2$ homology $\mathbb{H}_*(LS^2;\mathbb{F}_2)$ is not isomorphic
as Batalin-Vilkovisky algebras to the Hochschild cohomology of $H^*(S^2;\mathbb{F}_2)$,
$HH^{*}(H^*(S^2;\mathbb{F}_2);H^*(S^2;\mathbb{F}_2))$.
\end{cor}
This means exactly that there exists no isomorphism
between $\mathbb{H}_*(LS^2;\mathbb{F}_2)$ and
$HH^{*}(H^*(S^2;\mathbb{F}_2);H^*(S^2;\mathbb{F}_2))$
which at the same time,
\begin{itemize}
\item is an isomorphism of algebras and
\item commutes with the $\BV$ operators,
\end{itemize}
although separately
\begin{itemize}
\item there exists an isomorphism of algebras between $\mathbb{H}_*(LS^2;\mathbb{F}_2)$
and $HH^{*}(H^*(S^2;\mathbb{F}_2);H^*(S^2;\mathbb{F}_2))$
(Corollary~\ref{homologie des lacets egale a cohomologie de Hochschild
  comme G-algebre}) and
\item there exists also an isomorphism commuting with the $\BV$ operators
  between them.
\end{itemize}
\begin{proof}
By Proposition~\ref{Cohomologie de Hochschild des spheres modulo 2}, $HH^{*}(H^*(S^2);H^*(S^2))$ is
the Batalin-Vilkovisky algebra given by $\varepsilon=0$ in Lemma~\ref{BV algebre des lacets sur la sphere a epsilon pres modulo 2}.
On the contrary, by Theorem~\ref{homologie des lacets sur la sphere modulo 2},  $\mathbb{H}_*(LS^2;\mathbb{F}_2)$ is
the Batalin-Vilkovisky algebra given by $\varepsilon=1$.
At the end of the proof of Lemma~\ref{BV algebre des lacets sur la sphere a epsilon pres modulo 2},
we saw that the two cases $\varepsilon=0$ and $\varepsilon=1$ correspond to two non isomorphic
Batalin-Vilkovisky algebras.
\end{proof}
More generally, we believe that for any prime $p$, the free loop space modulo $p$ of the complex projective space
$\mathbb{H}_*(L\mathbb{CP}^{p-1};\mathbb{F}_p)$\footnote{B\"okstedt and
Ottosen~\cite{ConfDanemark}
have recently announced the computation of Batalin-Vilkovisky algebra
$\mathbb{H}_*(L\mathbb{CP}^{n};\mathbb{F}_p)$.}
is not isomorphic as Batalin-Vilkovisky algebras to the Hochschild
cohomology $HH^{*}(H^*(\mathbb{CP}^{p-1};\mathbb{F}_p);H^*(\mathbb{CP}^{p-1};\mathbb{F}_p))$.
Such phenomena for formal manifolds should not appear over a field of characteric $0$.

Recall that by Poincar\'e duality, we have the isomorphism
$$(\ref{dualite de Poincare})\quad\Theta:H^*(S^2)\buildrel{\cong}\over\rightarrow  H^*(S^2)^\vee.$$
Therefore we have the isomorphism
$$
HH^*(H^*(S^2);\Theta):HH^*(H^*(S^2);H^*(S^2))
\buildrel{\cong}\over\rightarrow  
HH^*(H^*(S^2);H^*(S^2)^\vee).
$$
Consider any isomorphism of graded algebras
\begin{equation}\label{Cohen-Jones}
\mathbb{H}_*(LS^2)\cong
HH^*(S^*(S^2);S^*(S^2)).
\end{equation}
By Corollary~\ref{homologie des lacets egale a cohomologie de Hochschild comme G-algebre}, such isomorphism exists.
Cohen and Jones (\cite[Theorem 3]{CohJon:ahtrostringtopology}
and~\cite{Cohen:multpropAtiyahduality}) proved that
such isomorphism exists for any manifold $M$.
Since $S^2$ is formal, we have the isomorphism of algebras
$$(\ref{invariance cohomologie de Hochschild})
\quad HH^*(S^*(S^2);S^*(S^2))\buildrel{\cong}\over\rightarrow
HH^*(H^*(S^2);H^*(S^2)).$$ 
By~\cite{JonesJ:Cycheh}, we
have the isomorphisms of $H_*(S^1)$-modules
$$
H_*(LS^2)\buildrel{(\ref{Jones})}\over\cong HH^*(S^*(S^2);S^*(S^2)^\vee)
\buildrel{(\ref{invariance homologie de Hochschild})}\over\cong HH^*(H^*(S^2);H^*(S^2)^\vee).
$$
Corollary~\ref{homologie des lacets different de cohomologie de Hochschild comme BV-algebre} implies that the following diagram does not commute over
$\mathbb{F}_2$:

\xymatrix{
&HH^*(S^*(S^2);S^*(S^2)^\vee)
\ar[r]^{(\ref{invariance homologie de Hochschild})}
&HH^*(H^*(S^2);H^*(S^2)^\vee)\\
H_*(LS^2)\ar[ur]^{(\ref{Jones})}
\ar[dr]^{(\ref{Cohen-Jones})}\\
& HH^*(S^*(S^2);S^*(S^2))
\ar[r]^{(\ref{invariance cohomologie de Hochschild})}
&HH^*(H^*(S^2);H^*(S^2))\ar[uu]_{HH^*(H^*(S^2);\Theta)}
}
This is surprising because as explained by Cohen and
Jones~\cite[p. 792]{CohJon:ahtrostringtopology},
the composite 
of the isomorphism (\ref{Jones}) given by
Jones in~\cite{JonesJ:Cycheh}
and an isomorphism induced by Poincar\'e duality
should give
an isomorphism of algebras between $\mathbb{H}_*(LS^2)$ and
$HH^*(S^*(S^2);S^*(S^2))$.
\section{Appendix by Gerald Gaudens and Luc Menichi.}
Let $X$ be a pointed topological space.
Recall that the circle $S^1$ acts on the double loop space
$\Omega^2 X$.
Consider the induced operator
$\BV:H_*(\Omega^2 X)\rightarrow H_{*+1}(\Omega^2 X)$.
Getzler~\cite{Getzler:BVAlg} has showed that $H_*(\Omega^2 X)$ equipped with the Pontryagin product and this operator $\BV$
forms a Batalin-Vilkovisky algebra.
In~\cite{Gaudens-Menichi}, Gerald Gaudens and the author have determined this Batalin-Vilkovisky algebra
$H_*(\Omega^2 S^3;\mathbb{F}_2)$.
The key was the following Theorem.
In~\cite[Proposition 7.46]{Kallel:livre}, answering to a question of
Gerald Gaudens,
Sadok Kallel and Paolo Salvatore give another proof of this Theorem.
\begin{theor}~\cite{Gaudens-Menichi}\label{Gaudens Menichi}
The operator $\BV:H_1(\Omega^2 S^3;\mathbb{F}_2)\rightarrow H_2(\Omega^2 S^3;\mathbb{F}_2)$ is non trivial.
\end{theor}
Both proofs~\cite{Gaudens-Menichi} and~\cite[Proposition 7.46]{Kallel:livre}
are unpublished and publicly unavailable yet.
So the goal of this section is to give a proof of this theorem which is
as simple as possible.

Denote by $*$ the Pontryagin product in $H_*(\Omega^2 X)$ and
by $\circ$ the map induced in homology by the composition map
$\Omega^2 X\times \Omega^2 S^2\rightarrow\Omega^2 X$.
Denote by $\Omega^2_n S^2$, the path-connected component of the
degree $n$ maps.
Denote by $v_1$ the generator of
$H_1(\Omega^2_0 S^2;\mathbb{F}_2)$ and by $[1]$ the generator of
$H_0(\Omega^2_1 S^2;\mathbb{F}_2)$. 
\begin{lem}\label{BV lacets double en general}
For $x\in H_*(\Omega^2 X;\mathbb{F}_2)$,
$
\BV x=x\circ (v_1*[1]).
$
\end{lem}
\begin{proof}
The circle $S^1$ acts on the sphere $S^2$.
Therefore we have a morphism of topological monoids
$
\Theta:(S^1,1)\rightarrow (\Omega^2_1 S^2,id_{S^2})
$.
The action of $S^1$ on $\Omega^2 X$ is the composite
$
S^1\times\Omega^2 X\buildrel{\Theta\times\Omega^2 X}\over\rightarrow
\Omega^2_1 S^2\times\Omega^2 X
\buildrel{\circ}\over\rightarrow\Omega^2 X
$.
Therefore for $x\in H_*(\Omega^2 X;\mathbb{F}_2)$,
$
\BV x=x\circ (H_1(\Theta)[S^1])
$.

Suppose that $H_1(\Theta)[S^1]=0$.
Then for any topological space $X$, the operator $\BV$
on $H_*(\Omega^2 X;\mathbb{F}_2)$ is null.
Therefore, for any $x$ and $y\in H_*(\Omega^2 X;\mathbb{F}_2)$,
$\{x,y\}=\BV(xy)-(\BV x)y-x(\BV y)=0$.
That is the modulo $2$ Browder brackets on any double loop space
are null. This is obviously false.
For example, Cohen in~\cite{CohenF:configshLiea} explains
that the Gerstenhaber algebra $H_*(\Omega^2 \Sigma^2 Y)$
has in general many non trivial Browder brackets.
So the assumption $H_1(\Theta)[S^1]=0$ is false.

Since the loop multiplication by $id_{S^2}$
in the $H$-group $\Omega^2 S^2$, is a homotopy equivalence,
the Pontryagin product by $[1]$,
$*[1]:H_*(\Omega^2_0 S^2)\buildrel{\cong}\over\rightarrow
H_*(\Omega^2_1 S^2)$ is an isomorphism.
Therefore $v_1*[1]$ is a generator of
$
H_1(\Omega^2_1 S^2)
$
So $H_1(\Theta)[S^1]=v_1*[1]$.
So finally
$$
\BV x=x\circ (H_1(\Theta)[S^1])=x\circ (v_1*[1]).
$$
\end{proof}
Recall that $v_1$ denote the generator of
$H_1(\Omega^2_0 S^2;\mathbb{F}_2)$.
\begin{lem}\label{Delta dans lacets doubles sur la sphere}
In the Batalin-Vilkovisky algebra $H_*(\Omega^2 S^2;\mathbb{F}_2)$,
$\BV(v_1)=v_1*v_1$.
\end{lem}
\begin{proof}
Recall that $[1]$ is the generator of $H_0(\Omega^2_1 S^2)$.
By Lemma~\ref{BV lacets double en general},
$$
\BV [1]=[1]\circ (v_1*[1])=(v_1*[1]).
$$
Denote by $Q:H_q(\Omega^2_n S^2)\rightarrow H_{2q+1}(\Omega^2_{2n} S^2)$
the Dyer-Lashof operation.
It is well known that $Q[1]=v_1*[2]$.
So by~\cite[Theorem 1.3 (4) p. 218]{Cohen-Lada-May:homiterloopspaces}
$$\{v_1*[2],[1]\}=\{Q[1],[1]\}=\{[1],\{[1],[1]\}\}.$$
By~\cite[Theorem 1.2 (3) p. 215]{Cohen-Lada-May:homiterloopspaces},
 $\{[1],[1]\}=0$.
Therefore on one hand, $\{v_1*[2],[1]\}$ is null.
And on the other hand, using the Poisson
relation~(\ref{Poisson relation}), since
$
\{[2],[1]\}=\{[1]*[1],[1]\}=2\{[1],[1]\}*[1]=0$,
$$
\{v_1*[2],[1]\}=\{v_1,[1]\}*[2]+v_1*\{[2],[1]\}=\{v_1,[1]\}*[2].
$$
Since $*[1]:H_*(\Omega^2 S^2)\buildrel{\cong}\over
\rightarrow H_*(\Omega^2 S^2)$
is an isomorphism, we obtain that Browder bracket $\{v_1,[1]\}$ is null.
Therefore,
$$\BV(v_1*[1])=(\BV v_1)*[1]+v_1*(\BV[1])=((\BV v_1)-v_1*v_1)*[1].$$
But $\BV(v_1*[1])=(\BV\circ\BV)([1])=0$. Therefore
$(\BV v_1)$ must be equal to $v_1*v_1$.
\end{proof}
\begin{proof}[Proof of Theorem~\ref{Gaudens Menichi}]
We remark that since $\BV$ preserves path-connected components
and since the loop multiplication of two homotopically trivial loops
is a homotopically trivial loop, $H_*(\Omega^2_0 S^2)$ is a sub
Batalin-Vilkovisky algebra of $H_*(\Omega^2 S^2)$.

Let $S^1\hookrightarrow S^3\buildrel{\eta}\over\twoheadrightarrow S^2$
be the Hopf fibration. After double looping, the Hopf fibration gives
the fibration
$\Omega^2 S^1\hookrightarrow \Omega^2 S^3
\buildrel{\Omega^2 \eta}\over\twoheadrightarrow \Omega^2_0 S^2$
with contractile fiber $\Omega^2 S^1$ and path-connected base
$\Omega^2_0 S^2$. Therefore $\Omega^2 \eta:\Omega^2 S^3
\buildrel{\simeq}\over\rightarrow \Omega^2_0 S^2$
is a homotopy equivalence.
And so $H_*(\Omega^2\eta):H_*(\Omega^2 S^3)
\buildrel{\cong}\over\rightarrow H_*(\Omega^2_0 S^2)$
is an isomorphism of Batalin-Vilkovisky algebras.

Let $u_1$ be the generator of $H_1(\Omega^2 S^3)$.
Lemma~\ref{Delta dans lacets doubles sur la sphere}
implies that $\BV(u_1)=u_1*u_1$.
Since $u_1*u_1$ is non zero in $H_*(\Omega^2 S^3;\mathbb{F}_2)$,
$\BV(u_1)$ is non trivial.
\end{proof}
\bibliography{Bibliographie}
\bibliographystyle{amsplain}
\end{document}